\documentclass[a4paper,11pt]{article}
\usepackage[english]{babel}
\usepackage{amsfonts}
\usepackage{amsthm}
\usepackage{amsbsy}
\usepackage{amssymb}
\usepackage{graphicx}
\usepackage{eso-pic}
\usepackage{dsfont}
\usepackage{tikz}
\usepackage{xfrac}
\usepackage{float}
\usepackage{bookmark}
\usepackage{amsmath}
\usepackage{hyperref}
\usepackage{lscape}
\usetikzlibrary{matrix,arrows}
\usepackage[letterpaper,left=0.8in,right=0.8in,top=.75in,bottom=.75in]{geometry}
\newtheorem{lemma}{Lemma}[section]
\newtheorem{teo}[lemma]{Theorem}
\newtheorem{prop}[lemma]{Proposition}
 
\newtheorem{conj}[lemma]{Conjecture}

\theoremstyle{definition}

\newtheorem{quest}[lemma]{Question}

\theoremstyle{remark}

\begin{document}
\title{Spherical tetrahedra with rational volume, \\  and spherical Pythagorean triples}
\author{Alexander Kolpakov \& Sinai Robins}
\date{}
\maketitle

\begin{abstract}
\noindent   We study spherical tetrahedra with rational dihedral angles and rational volumes.  Such tetrahedra occur in the Rational Simplex Conjecture by Cheeger and Simons, and we supply vast families, discovered by computational  efforts, of positive examples that confirm this conjecture. As a by-product, we also obtain a classification of all spherical Pythagorean triples, previously found by Smith. 

\end{abstract}

\section{Introduction}\label{section:pre}

A \textit{spherical tetrahedron} $T$ can be defined as the intersection of a simplicial cone in $\mathbb R^4$ with the unit sphere $\mathbb{S}^3$ centred at the origin.   In other words, $T$ has four vertices connected by spherical geodesics on $\mathbb{S}^3$ that comprise its edges, and each of its vertices is the intersection
of exactly three of its spherical facets.  A \textit{spherical Coxeter tetrahedron} $T$ is a spherical tetrahedron whose six dihedral angles are of the form $\pi/n$, with $n\geq 2$.

A complete list of spherical Coxeter tetrahedra was produced by Coxeter \cite{Coxeter}, and shows that there are eleven types of spherical Coxeter tetrahedra in $\mathbb{S}^3$. Let $S_i$, $i=1,\dots,11$, denote these spherical tetrahedra, as presented in Table~\ref{tabular:coxeter-tetrahedra}. 

In the present paper we study \textit{rational spherical tetrahedra}, as generalisations of spherical Coxeter tetrahedra, where we now allow their dihedral angles to be arbitrary rational multiples of $\pi$.  An important focus here is the determination of their volume, which is also called a solid angle in some of the literature.
 
The volume of a spherical Coxeter tetrahedron is easily seen to be a rational multiple of the total volume of the sphere $\mathbb{S}^3$, which is $2 \pi^2$. We describe a wide class of  rational spherical tetrahedra whose volumes are rational multiples of $\pi^2$, in relation to the work of Cheeger and Simons \cite{CS}. 

\begin{table}[h]
\begin{centering}
\begin{tabular}{c|c|c|c}
$i$& Symbol& Coxeter diagram& Volume\\
\hline
1& $A_4$& 
	\begin{tikzpicture}
	\draw (0,0) -- (1,0);
    \draw (1,0) -- (2,0);
	\draw (2,0) -- (3,0);
	
	\draw[fill=white] (0,0) circle(.1);
	\draw[fill=white] (1,0) circle(.1);
	\draw[fill=white] (2,0) circle(.1);
	\draw[fill=white] (3,0) circle(.1);
	\end{tikzpicture}&
	$\frac{\pi^2}{60}$\\
\hline
2& $B_4$& 
	\begin{tikzpicture}
	\draw (0,0) -- (1,0);
    \draw (1,0) -- (2,0);
	\draw (2,0.05) -- (3,0.05);
	\draw (2,-0.05) -- (3,-0.05);
	
	\draw[fill=white] (0,0) circle(.1);
	\draw[fill=white] (1,0) circle(.1);
	\draw[fill=white] (2,0) circle(.1);
	\draw[fill=white] (3,0) circle(.1);
	\end{tikzpicture}&
	$\frac{\pi^2}{192}$\\
\hline
3& $D_4$& 
	\begin{tikzpicture}
	\draw (0,0) -- (1,0);
	\draw (1,0) -- (2,0.5);
	\draw (1,0) -- (2,-0.5);
	
	\draw[fill=white] (0,0) circle(.1);
	\draw[fill=white] (1,0) circle(.1);
	\draw[fill=white] (2,0.5) circle(.1);
	\draw[fill=white] (2,-0.5) circle(.1);
	\end{tikzpicture}&
	$\frac{\pi^2}{96}$\\
\hline
4& $H_4$& 
	\begin{tikzpicture}
	\draw (0,0) -- (1,0);
    \draw (1,0) -- (2,0);
	\draw (2,0.05) -- (3,0.05);
	\draw (2,0) -- (3,0);
	\draw (2,-0.05) -- (3,-0.05);
	
	\draw[fill=white] (0,0) circle(.1);
	\draw[fill=white] (1,0) circle(.1);
	\draw[fill=white] (2,0) circle(.1);
	\draw[fill=white] (3,0) circle(.1);
	\end{tikzpicture}&
	$\frac{\pi^2}{7200}$\\
\hline
5& $F_4$& 
	\begin{tikzpicture}
	\draw (0,0) -- (1,0);
    \draw (1,0.05) -- (2,0.05);
	\draw (1,-0.05) -- (2,-0.05);
	\draw (2,0) -- (3,0);
	
	\draw[fill=white] (0,0) circle(.1);
	\draw[fill=white] (1,0) circle(.1);
	\draw[fill=white] (2,0) circle(.1);
	\draw[fill=white] (3,0) circle(.1);
	\end{tikzpicture}&
	$\frac{\pi^2}{576}$\\
\hline
6& $A_3\times A_1$& 
	\begin{tikzpicture}
	\draw (0,0) -- (1,0);
    \draw (1,0) -- (2,0);
	
	\draw[fill=white] (0,0) circle(.1);
	\draw[fill=white] (1,0) circle(.1);
	\draw[fill=white] (2,0) circle(.1);
	\draw[fill=white] (3,0) circle(.1);
	\end{tikzpicture}&
	$\frac{\pi^2}{24}$\\
\hline
7& $B_3\times A_1$& 
	\begin{tikzpicture}
	\draw (0,0) -- (1,0);
	\draw (1,0.05) -- (2,0.05);
    \draw (1,-0.05) -- (2,-0.05);
	
	\draw[fill=white] (0,0) circle(.1);
	\draw[fill=white] (1,0) circle(.1);
	\draw[fill=white] (2,0) circle(.1);
	\draw[fill=white] (3,0) circle(.1);
	\end{tikzpicture}&
	$\frac{\pi^2}{48}$\\
\hline
8& $H_3\times A_1$& 
	\begin{tikzpicture}
	\draw (0,0) -- (1,0);
	\draw (1,0.05) -- (2,0.05);
	\draw (1,0) -- (2,0);
    \draw (1,-0.05) -- (2,-0.05);
	
	\draw[fill=white] (0,0) circle(.1);
	\draw[fill=white] (1,0) circle(.1);
	\draw[fill=white] (2,0) circle(.1);
	\draw[fill=white] (3,0) circle(.1);
	\end{tikzpicture}&
	$\frac{\pi^2}{120}$\\
\hline
9& $I_2(k)\times I_2(l)$& 
	\begin{tikzpicture}
	\draw (0,0) -- (1,0);
	\draw (2,0) -- (3,0);
	
	\draw[fill=white] (0,0) circle(.1);
	\draw[fill=white] (1,0) circle(.1);
	\draw[fill=white] (2,0) circle(.1);
	\draw[fill=white] (3,0) circle(.1);
	
	\node at (0.5, 0.3) {$k$};
	\node at (2.5, 0.3) {$l$};
	\end{tikzpicture}&
	$\frac{\pi^2}{2kl}$\\
\hline
10& $I_2(k)\times A_1^{\times 2}$& 
	\begin{tikzpicture}
	\draw (0,0) -- (1,0);
	
	\draw[fill=white] (0,0) circle(.1);
	\draw[fill=white] (1,0) circle(.1);
	\draw[fill=white] (2,0) circle(.1);
	\draw[fill=white] (3,0) circle(.1);
	
	\node at (0.5, 0.3) {$k$};
	\end{tikzpicture}&
	$\frac{\pi^2}{4k}$\\
\hline
11& $A_1^{\times 4}$& 
	\begin{tikzpicture}	
	\draw[fill=white] (0,0) circle(.1);
	\draw[fill=white] (1,0) circle(.1);
	\draw[fill=white] (2,0) circle(.1);
	\draw[fill=white] (3,0) circle(.1);
	\end{tikzpicture}&
	$\frac{\pi^2}{8}$\\
\end{tabular}
\caption{Coxeter tetrahedra in $\mathbb{S}^3$}
\label{tabular:coxeter-tetrahedra}
\end{centering}
\end{table}

In this work, an angle $\alpha$ (assumed to be a plane angle of a polygon, or a dihedral angle of a polyhedron) is called \textit{rational} if $\alpha \in \pi\,\mathbb{Q}$. 
Similarly, an edge of a polygon (or an edge length of a polyhedron) of length $l$ is called \textit{rational}
 if $l \in \pi \, \mathbb{Q}$.  
 Finally, an $n$-tuple of numbers $(x_1, \dots, x_n)$ is \textit{rational}, if $x_i \in \pi\,\mathbb{Q}$ for every $1 \leq i \leq n$. 

Descending to      $\mathbb{S}^2 \subset \mathbb R^3$,
a {\em spherical Pythagorean triple}  is defined to be a rational solution $(p, q, r)$ to 
\begin{equation}\label{eq:Pythagorean-triple}
\cos p \cdot \cos q + \cos r = 0,
\end{equation}
where $\pi - p$, $\pi - q$, and $\pi - r$ are the side lengths of a spherical right triangle $T$. 
The side lengths of a spherical triangle are subject to several additional constraints on $p$, $q$ and $r$:
\begin{equation*}
0<\, p,\,\, q,\,\, r\, < \pi, \hspace {0.25in} p + q + r < 2\pi,
\end{equation*}
\begin{equation*}
p + q < r,\hspace {0.5in} p + r < q,\hspace {0.5in} q + r < p.
\end{equation*}
We relax the above conditions and call any solution of \eqref{eq:Pythagorean-triple},  with $0 < p, q, r < \pi$, a Pythagorean triple. 

\begin{quest}\label{q1}
Is there any reasonably simple classification of rational Pythagorean triples, corresponding to the side lengths of a spherical triangle? 
\end{quest}

\begin{figure}[ht]
\center
\includegraphics[scale=0.4]{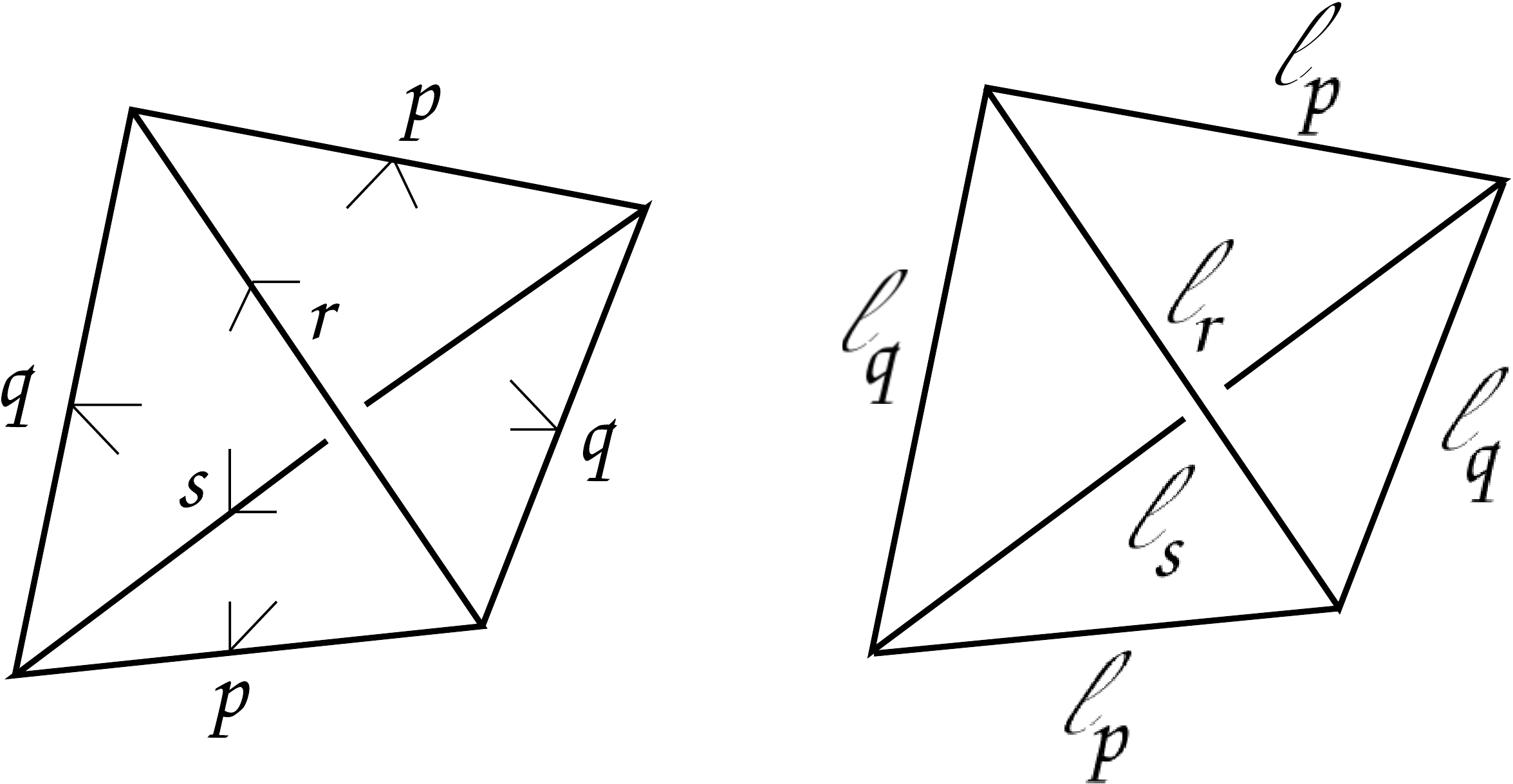}
\caption{The dihedral angles (left) and edge lengths (right) of a $\mathbb{Z}_2$-symmetric tetrahedron $T$.}\label{fig:z2-symmetric-tetrahedron}
\end{figure}

Returning to $\mathbb{S}^3$, we focus on a broader class of  ``Pythagorean quadruples'', that will become useful in the discussion of $\mathbb{Z}_2$-symmetric spherical tetrahedra with rational dihedral angles (or \textit{rational tetrahedra}, for short) later on. To this end, we call $(p, q, r, s)$ a \textit{spherical Pythagorean quadruple} if it is a solution to the equation
\begin{equation}\label{eq:Pythagorean-quadruple}
\cos p \cdot \cos q+ \cos \frac{r+s}{2}\cdot \cos \frac{r-s}{2} = 0.
\end{equation}
Here, we shall suppose that $0 < p, q, r, s < \pi$. The corresponding spherical tetrahedron, if it exists, looks akin to the one depicted in Figure~\ref{fig:z2-symmetric-tetrahedron} and is called a $\mathbb{Z}_2$-symmetric (spherical) tetrahedron. 

\noindent
We note that a quadruple with $r = s$ corresponds to the usual Pythagorean triple $(p, q, r)$.

\begin{quest}\label{q2}
Is there any reasonably simple classification of rational Pythagorean quadruples corresponding to the dihedral angles of a spherical tetrahedron? 
\end{quest}

We shall answer Questions \ref{q1} and \ref{q2} simultaneously by classifying all Pythagorean quadruples.

\begin{teo}\label{thm:Pythagorean}
There exist exactly $59$ sporadic Pythagorean quadruples, and $42$ continuous families of Pythagorean quadruples corresponding to the dihedral angles of a $\mathbb{Z}_2$-symmetric spherical tetrahedron. 
\end{teo}

The proof of Theorem~\ref{thm:Pythagorean} is contained in Section 2.1 for the case of sporadic instances listed in Appendix A, and in Section 2.2 for the case of continuous families listed in Appendix B. 
The main tool in our proof is a very basic enumeration realised by a \texttt{SageMath} script Monty \cite{Monty}. Thus, Theorem~\ref{thm:Pythagorean} extends a result (unpublished) of Smith \cite{Smith} that classifies rational spherical Pythagorean triples by using a beautiful geometric connection with \textit{three-dimensional} Coxeter simplices: 

\begin{teo}[Theorem 2 in \cite{Smith}] \label{thm:Smith}
Aside from the trivial continuous family of solutions $(\pi/2, b, \pi/2)$, $0\leq b \leq \pi/2$, there is exactly one solution $(a,b,c)$ to $\cos a \cos b = \cos c$ with $0 \leq a, b, c \leq \pi/2$ being rational multiples of $\pi$, namely $(\pi/4, \pi/4, \pi/3)$.
\end{teo}
In terms of equation \eqref{eq:Pythagorean-triple}, the non-trivial triple in the above theorem is $(\pi/4, \pi/4, 2\pi/3)$.  All rational Pythagorean triples found by Smith belong to the continuous families of quadruples described in Section \ref{section:spherical-triples}.

The following statement is an observation which had its origin in the list of rational spherical 
Pythagorean quadruples, and which is of interest in the context of \cite{Felikson1, Felikson2}.

\begin{teo}\label{thm:rational-not-coxeter-tetrahedron}
There exists a rational tetrahedron in $\mathbb{S}^3$ whose volume has a value in $\pi^2\,\mathbb{Q}$, and which is not decomposable into any finite number of spherical Coxeter tetrahedra. 
\end{teo}

Thus, we can show that the property of ``being rational'' for a spherical tetrahedron is very far from ``being Coxeter'', even if its volume is a rational multiple of $\pi^2$, which is always true for Coxeter tetrahedra in $\mathbb{S}^3$. Here we recall that $\mathbb{S}^3$ has volume $2\pi^2$ in its natural metric of constant sectional curvature $+1$, and that every Coxeter polyhedron in $\mathbb{S}^3$ is a tetrahedron, which generates a finite discrete reflection group by reflection in its faces. 

The first open problem that Theorem \ref{thm:rational-not-coxeter-tetrahedron} vaguely relates to is Schl\"{a}fli's Conjecture: 
\begin{conj}(Schl\"{a}fli)
Let $T$ be an orthoscheme in $\mathbb{S}^3$ with rational dihedral angles. Then the volume of $T$ takes values in $\pi^2\, \mathbb{Q}$ if and only if $T$ is a Coxeter orthoscheme. 
\end{conj} 
The above statement can be generalised for spherical simplices of dimension $\geq 4$, and this is how it actually appears in Scl\"afli's original work \cite[p.~267, Formeln (4)--(5)]{Schlafli}.
Here, an orthoscheme is a tetrahedron with three mutually orthogonal faces that do not share a common vertex.  However, the tetrahedron mentioned in Theorem \ref{thm:rational-not-coxeter-tetrahedron} is not an orthoscheme. 

Another related open problem is the following question posed in \cite{CS} by Cheeger and Simons, and known as the Rational Simplex Conjecture:
\begin{quest}
\textit{Is it true that the volume of a rational spherical tetrahedron always takes values in $\pi^2\,\mathbb{Q}$?}
\end{quest}
The putative answer would be negative for ``virtually all'' rational simplices. Our results only show that the Rational Simplex Conjecture may hold for a tetrahedron which is geometrically ``far enough'' from a Coxeter tetrahedron, and thus one may still expect many ``positive examples''. Finally, we can produce many pairs of non-isometric rational tetrahedra with equal volumes and Dehn invariants. In view of Hilbert's 3\textsuperscript{rd} problem, it would be natural to ask if our examples are scissors congruent. 

\begin{center}
\textbf{Acknowledgements}\\
\end{center}
\noindent
{\small The authors gratefully acknowledge the support they received from Brown University and the Institute for Computational and Experimental Research in Mathematics - ICERM (Providence, USA), where most of this work has been conceived and carried out, and especially thank the ICERM program ``Point Configurations in Geometry, Physics and Computer Science'' for hospitality and excellent working atmosphere during their stay. A.K. has been supported by the Swiss National Science Foundation - SNSF project no.~PP00P2-170560, S.R. has been supported by the S\~{a}o Paulo Research Foundation - FAPESP project no.~15/10323-7. The authors are grateful to Don Zagier (MPIM, Germany and ICTP, Italy), John Parker (Durham University, UK), Richard Schwartz (Brown University, USA) and Ruth Kellerhals (Universit\'e de Fribourg, Switzerland) for stimulating discussions. They also thank the anonymous referees for their valuable comments and suggestions. }


\section{Pythagorean quadruples}\label{section:spherical-triples}

Let a spherical tetrahedron $T$ be defined as an intersection of a simplicial cone $C$ in $\mathbb{R}^4$ centred at the origin with the unit sphere $\mathbb{S}^3 = \{  \mathbf{v} = (x_1, x_2, x_3, x_4) \in \mathbb{R}^4\, |\, \|  \mathbf{v} \| = 1 \}$. We suppose that the dihedral angles of $C$ belong to the interval $(0, \pi)$. 

The dihedral angles of $T$ are equal to the corresponding dihedral angles between the three-dimensional faces of its defining cone $C$ measured at its two-dimensional faces. The edge lengths of $T$ correspond to the plane angles in the two-dimensional proper sub-cones  measured at the origin.

The polar dual $T^*$ of a spherical tetrahedron $T$, defined by a cone $C$, is the intersection of the dual cone $C^*$ with $\mathbb{S}^3$. 

We recall that a spherical tetrahedron is called $\mathbb{Z}_2$-symmetric, if it admits such a distribution of dihedral angles values as shown in Figure~\ref{fig:z2-symmetric-tetrahedron}. A Pythagorean quadruple of dihedral angles $(p, q, r, s)$ of a $\mathbb{Z}_2$-symmetric spherical tetrahedron is a solution to equation \eqref{eq:Pythagorean-quadruple}.

Then, by cosidering polar duals, one can deduce from Proposition~6 of \cite{KMP} the following:

\begin{prop}\label{prop:volume-angles}
If $p$, $q$, $r$ and $s$ are the dihedral angles of a $\mathbb{Z}_2$-symmetric spherical tetrahedron $T$, for which equation \eqref{eq:Pythagorean-quadruple} holds, then the volume of $T$ can be expressed as
\begin{equation}\label{eq:volume-angles}
\mathrm{Vol}\, T = \frac{1}{2}\left( \frac{r\,(2\pi-r)}{2} + p^2 + q^2 + \frac{s\,(2\pi-s)}{2} - \pi^2  \right).
\end{equation}
\end{prop}

Thus, once the dihedral angles of a tetrahedron $T$ as above are rational, then its volume has a value in $\pi^2\,\mathbb{Q}$. It also follows from \cite[Proposition 6]{KMP} (and the discussion preceding it), that a rational $\mathbb{Z}_2$-symmetric tetrahedron has rational edge lengths. Namely, the following holds.

\begin{prop}\label{prop:edge-lengths}
If $(p, q, r, s)$ is the quadruple of dihedral angles of a $\mathbb{Z}_2$-symmetric spherical tetrahedron $T$, for which equation \eqref{eq:Pythagorean-quadruple} holds, then the lengths of its respective edges, as depicted in Figure~\ref{fig:z2-symmetric-tetrahedron} are given by the quadruple $(\ell_p, \ell_q, \ell_r, \ell_s) = (p, q, \pi - r, \pi - s)$. 
\end{prop}

Once we have $r = s$ for a spherical $\mathbb{Z}_2$-symmetric tetrahedron $T$, we get a triple $(p, q, r)$, which corresponds in this case to a symmetric spherical tetrahedron, rather than to a triangle. However, $(p,q,r)$ is a Pythagorean triple in the sense of our initial definition. Indeed, for each vertex $v$ of $T$ in this case, its link $\mathrm{Lk}_v$ is a rational spherical triangle with plane angles $p$, $q$, and $r$. Its dual $\mathrm{Lk}^*_v$ is a spherical triangle with edge length $\pi - p$, $\pi - q$, $\pi - r$, while $p$, $q$, and $r$ satisfy equation \eqref{eq:Pythagorean-triple}.

A Pythagorean quadruple $(p, q, r, s)$ represents the dihedral angles of a $\mathbb{Z}_2$-symmetric spherical tet\-ra\-hed\-ron $T$, if and only if the associated Gram matrix 
\begin{equation}\label{eq:Gram-matrix}
G = G(T) := \left( \begin{array}{cccc}
1& -\cos r& -\cos p& -\cos q\\
-\cos r& 1& -\cos q& -\cos p\\
-\cos p& -\cos q& 1& -\cos s\\
-\cos q& -\cos p& -\cos s& 1
\end{array}  \right)
\end{equation}
is positive definite \cite[Lemma~1.2]{Luo}. 

Thus, once we have a rational solution $(p, q, r, s)$ to \eqref{eq:Pythagorean-quadruple}, then we only need to check if the Gram matrix $G(T)$ given by \eqref{eq:Gram-matrix} is positive definite. If it is indeed the case, then we obtain a rational spherical tetrahedron $T$ such that $\mathrm{Vol}\,T \in \pi^2\, \mathbb{Q}$. 

First of all, finding a solution to equation \eqref{eq:Pythagorean-quadruple} is equivalent to finding a solution to the equation
\begin{equation}\label{eq:quad-angles-1}
\cos(a) + \cos(b) + \cos(c) + \cos(d) = 0,
\end{equation}  
while the correspondence between two sets of solutions is given by
\begin{equation}\label{eq:quad-angles-2}
p = \frac{a+b}{2}, \,\, q = \frac{a-b}{2}, \,\, r = c,\,\, s = d.
\end{equation}

We shall search for all possible solutions to \eqref{eq:quad-angles-1}  -- \eqref{eq:quad-angles-2}, such that $0 < p, q, r, s < \pi$, and $r \geq s$. The former condition is necessary for the dihedral angles of a spherical tetrahedron $T$, and the latter can be assumed since $r$ and $s$, as drawn in Figure~\ref{fig:z2-symmetric-tetrahedron}, can be interchanged by an isometry of $\mathbb{S}^3$ without interchanging $p$ and $q$.

If $(a, b, c, d)$ is a rational quadruple, then \eqref{eq:quad-angles-1} turns out to be a trigonometric Dio\-phan\-tine equation which has been studied by Conway and Jones in \cite{CJ}. All of its solutions such that $0 < a$, $b$, $c$, $d < \frac{\pi}{2}$ are listed in \cite[Theorem~7]{CJ}. 
For convenience, its statement is reproduced below, although using a slightly different notation.
\begin{teo}[Theorem 7 in \cite{CJ}]\label{thm:cj}
Suppose that we have at most four rational multiples of $\pi$ lying strictly between $0$ and $\pi/2$ for which some rational linear combination $S$ of their cosines is rational, but no proper subset has this property. Then $S$ is proportional to one of the following list:
\begin{itemize}
\item[1.] $\cos \frac{\pi}{3} - \cos \frac{\pi}{3}$ ($=0$),
\item[2.] $-\cos t + \cos \left(  t + \frac{\pi}{3} \right) + \cos \left(  t - \frac{\pi}{3} \right)$ ($=0$),
\item[3.] $\cos \frac{\pi}{5} - \cos \frac{2\pi}{5} - \cos \frac{\pi}{3}$ ($=0$),
\item[4.] $\cos \frac{\pi}{7} - \cos \frac{2\pi}{7} + \cos \frac{3\pi}{7} - \cos \frac{\pi}{3}$ ($=0$),
\item[5.] $ \cos \frac{\pi}{5} - \cos \frac{\pi}{15} + \cos \frac{4\pi}{15} - \cos \frac{\pi}{3}$ ($=0$),
\item[6.] $-\cos \frac{2\pi}{5} + \cos \frac{2\pi}{15} - \cos \frac{7\pi}{15} - \cos \frac{\pi}{3}$ ($=0$),
\item[7.] $\cos \frac{\pi}{7} + \cos \frac{3\pi}{7} - \cos \frac{\pi}{21} + \cos \frac{8\pi}{21}$ $\left( =\frac{1}{2} \right)$,
\item[8.] $\cos \frac{\pi}{7} - \cos \frac{2\pi}{7} + \cos\frac{2\pi}{21} - \cos \frac{5\pi}{21}$ $\left( =\frac{1}{2} \right)$,
\item[9.] $- \cos \frac{2\pi}{7} + \cos \frac{3\pi}{7} + \cos \frac{4\pi}{21} + \cos \frac{10 \pi}{21}$ $\left( = \frac{1}{2} \right)$,
\item[10.] $- \cos \frac{\pi}{15} + \cos \frac{2\pi}{15} + \cos \frac{4\pi}{15} - \cos \frac{7\pi}{15}$ $\left( = \frac{1}{2} \right)$.
\end{itemize}
\end{teo}

According to Theorem~\ref{thm:cj} there is a single continuous family of linear combinations of cosines, depending on a real-valued parameter $t$, which, for every instance of $t \in \pi \mathbb{Q}$, provides a rational solution to \eqref{eq:quad-angles-1}. The remaining linear combinations we call sporadic, in order to distinguish them from continuous families. Also, our methods to handle sporadic solutions to \eqref{eq:quad-angles-1} and their continuous families will be slightly different, since the former require more computations to be performed (first, numerically, and then exactly by verifying the respective minimal polynomials), while the latter need more symbolic algebra and the use of \texttt{SymPy} \cite{SymPy}.  

\subsection{Rational spherical tetrahedra: $59$ sporadic instances}\label{section:sporadic}

Let the rational length of a quadruple $(a, b, c, d)$ giving rise to the trigonometric sum $S = \cos a + \cos b + \cos c + \cos d$ in \eqref{eq:quad-angles-1} be defined as the maximal length of its sub-sum $S'$, such that $S' \in \mathbb{Q}$, but for any proper sub-sum $S''$ of $S'$ it still holds that $S'' \notin \mathbb{Q}$. 

Then,  we can already notice that there is no solution to \eqref{eq:quad-angles-1} of rational length $4$. Indeed, each linear combination of rational length $4$ would yield an expression $S$ equal to the right-hand side of items $7$, $8$, $9$, or $10$ in Theorem \ref{thm:cj}, up to a sign. None of those sums evaluates to $0$.

The sporadic solutions to \eqref{eq:quad-angles-1} mentioned in items $4$, $5$, and $6$ of Theorem \ref{thm:cj} have rational length $3$. The one mentioned in item 3 has rational length $2$. Finally, only those solutions where each cosine term of $S$ above is a rational number have rational length $1$. The latter is possible only if $a, b, c, d \in \{ 0, \pi/3, \pi/2 \}$, given that $0 \leq a$, $b$, $c$, $d \leq \frac{\pi}{2}$ . 

However, Theorem~\ref{thm:cj} provides only the sub-sums realising the rational length of $S$, and says nothing about the remaining part of the sum, which may have itself various rational length (e.g. if $S$ has rational length $2$ realised by a sub-sum $S'$, then $S - S'$ may have rational length $1$ or $2$). 

We shall need a wider range of dihedral angles represented by the Pythagorean quadruple $(a, b, c, d)$, namely $0 < a$, $b$, $c$, $d < \pi$. Thus, for each dihedral angle in each entry on the list of Theorem~\ref{thm:cj}, we also consider its complement to $\pi$ and $2\pi$, respectively. However, we always keep in mind that any angle in the interval $(0,\pi)$, as above, can be brought to an angle in $(0, \pi/2)$ in such a way that we do not create any new sums as compared to Theorem~\ref{thm:cj}, since the only difference will be some cosines in $S$ changing their signs. 

Moreover, if we assume that $a$, $b$, $c$, $d \in (0, \pi)$ instead of $(0, \pi/2)$, we need to consider one more continuous family in addition to the ones already mentioned in Theorem~\ref{thm:cj}. Namely, we need to consider $\cos \alpha + \cos \beta = 0$, with $\alpha = t$, $\beta = \pi - t$, and $t \in (0, \pi)$, as well as all possible complements of $\alpha$ and $\beta$ to $\pi$ and $2\pi$.

In order to simplify our search algorithm (at the cost of making it overall less efficient), we shall for each rational length of $S$ look at the set of possible denominators of the angles involved in $S'$ realising said length, and at the set of denominators realising any possible rational length of $S - S'$. Then we shall obtain a list of possible denominators $\delta_a$, $\delta_b$, $\delta_c$, $\delta_d$ that $a = \frac{\nu_a}{\delta_a} \pi$,  $b = \frac{\nu_b}{\delta_b} \pi$, $c = \frac{\nu_c}{\delta_c} \pi$, $d = \frac{\nu_d}{\delta_d} \pi$ may have, and choose their numerators $\nu_a$, $\nu_b$, $\nu_c$, $\nu_d$ so that $0 < a$, $b$, $c$, $d < \pi$. If any number of the form $\frac{\nu}{\delta} \pi$ equals $0$, then we assume $\delta = \infty$. Such an approach is still practically reasonable, and takes about $90$ minutes in total to run in \texttt{SageMath} \cite{Sage} on a MacPro 2.3 GHz Intel Core i5 Processor with 8 Mb RAM. 

An observation from Galois theory implies that if $S = \cos a + \cos b + \cos c + \cos d$ has rational length $1$, then the list of possible denominators of angles in $S$ is $L_0= \{ 1, 2, 3, \infty \}$.  

If $S$ has rational length $2$ realised by a sub-sum $S'$, then the list of possible denominators in $S'$ is $L_1 = \{ 3, 5 \}$ as indicated by item $3$ of Theorem \ref{thm:cj}, while the denominators in $S - S'$ can belong either to $L_0$ or to $L_1$.

If $S$ has rational length $3$ realised by a sub-sum $S'$, then the denominators of angles in $S'$ belong either to the list $L_2 = \{ 3, 7 \}$, or to $L_3 = \{ 3, 5, 15 \}$, as indicated by items $4$, $5$ and $6$ of Theorem \ref{thm:cj}, while the denominators of the remaining term $S - S'$ belong to $L_0$.

In Monty \cite{Monty} we use a brute-force search over the set of all dihedral angles with denominators from the union of the above mentioned lists $L_i$, $i \in \{0, 1, 2, 3\}$. This does not result in an \`a priori efficient search, however turns out to be sufficient to find all sporadic  solutions to \eqref{eq:quad-angles-1} and, subsequently, to \eqref{eq:Pythagorean-quadruple}. 

Each time a ``numerical'' zero is obtained in Monty's search, i.e. the condition $|S| < 10^{-8}$ is satisfied (which is a very generous margin for a numerical zero, since Monty's machine precision is $10^{-16}$), the minimal polynomial for $S$ is computed. Since $S$ is an algebraic integer, this test is sufficient to verify that $S = 0$. 

In each of the cases above, we check if the resulting dihedral angles $p$, $q$, $r$, $s$ of a ``candidate'' tetrahedron $T$ belong to the interior of the interval $(0, \pi)$, and whether the corresponding Gram matrix $G = G(T)$ is positive definite. The former condition guarantees that the first two corner minors $G_1 = 1$ and $G_2 = \sin^2 r$ of $G$, respectively of rank $1$ and $2$, are positive, and we need specifically to check only $G_3$ and $G_4 = \det G$. In Monty's search, $G_i$ is considered positive if $G_i > 10^{-8}$, which is again a generous numerical margin to decide if a number is positive. In order to verify that no possible solution is left out, we check if $G_i$ within the $10^{-8}$-neighbourhood of $0$ is actually $0$, by using minimal polynomials. Otherwise, $G_i < - 10^{-8}$, and is indeed negative. 

Finally, Monty finds $172$ sporadic solutions. Since the dihedral angles of all the listed tetrahedra satisfy equation \eqref{eq:Pythagorean-quadruple}, their volumes are rational multiples of $\pi^2$ by Proposition~\ref{prop:volume-angles}. 

There are, however, some of the sporadic solutions which belong by chance to one of the $42$ continuous families described in the next section. For brevity, we exclude them from our final list, and only $59$ genuinely sporadic solutions are presented in Appendix A. 

\subsection{Rational spherical tetrahedra: $42$ continuous families}\label{section:families}

By using a method analogous to the above, we find $34$ one-parameter continuous families, and $8$ two-parameter continuous families of rational spherical tetrahedra whose volumes take values in $\pi^2 \mathbb{Q}$. Those families are listed in Appendix B. When dealing with symbolic computations in Monty \cite{Monty}, we employ \texttt{SymPy} \cite{SymPy} in order to simplify expressions and check whether $S = 0$, rather than using the minimal polynomial test. 

In the case of continuous families, we have only two types of sub-sums $S'$ appearing in $S$, which depend on a parameter: 
\begin{itemize}
\item[i.] either a sub-sum of the form indicated in item 2 of Theorem \ref{thm:cj}, 
\item[ii.] or a sub-sum of the form $S'(t) = \cos(t) - \cos(t) = \cos(t) + \cos(\pi - t)$.
\end{itemize}

In the former case three of the angles $a$, $b$, $c$ and $d$ is \eqref{eq:quad-angles-1} belong to the list $L_0 = \{ \pi/3 - t, \pi/3 + t, 2\pi/3 - t, 2\pi/3 + t, \pi - t, t, \pi + t, 5 \pi/3 - t, 5\pi/3 + t \}$, with $t \in (0, \pi/6)$, and the remaining one belongs to $L_1 = \{ \pi/2, 3\pi/2 \}$. In the latter case, one pair of angles from $a$, $b$, $c$ and $d$ equals $\{ t, \pi - t \}$, with $t \in (0, \pi)$, and the remaining pair equals $\{ s, \pi - s  \}$, with $s \in (0, \pi)$. 

In case (i), we choose to produce graphs of the minors $G_3$ and $G_4 = \det G$ of the Gram matrix $G$ of each candidate tetrahedron, in order to check their positivity. The ones that appear positive on the whole interval $(0, \pi/6)$ indeed turn to $0$ only at the ends, or only one of the ends of the interval $(0, \pi/6)$. Then we check that those which appear negative on the interval $(0, \pi/6)$ do not turn positive near the end-points $0$ and $\pi/6$, but at worst become equal to $0$ at one or both of them. In order to verify all the above mentioned inequalities we use interval arithmetic implemented in \texttt{SageMath}, and for the equalities we use minimal polynomials, as before.  

In case (ii), we know that the tetrahedron $T^\ast$ with Coxeter diagram $A^{\times 4}_1$ belongs to any possible continuous family. The tetrahedron $T^\ast$ has all right angles, and thus  the minors $G_3(\pi/2, \pi/2)$ and $G_4(\pi/2, \pi/2)$ have to be positive for any family containing geometrically realisable tetrahedra. This filter leaves us with only few possible families, for which $G_3(s, t)$ and $G_4(s, t)$ have very simple form, amenable to elementary analysis for determining their positivity domains.

Finally, case (i) produces $34$ continuous families of tetrahedra depending on a single parameter, and case (ii) produces $8$ continuous families of tetrahedra depending on two parameters. All of them are listed in Appendix B, together with the domains of admissible parameter values, and the corresponding volume formulas. 


\section{Splitting rational polytopes into Coxeter tetrahedra}

Below we give a proof of Theorem~\ref{thm:rational-not-coxeter-tetrahedron}.  We begin by considering more closely one of the many Pythagorean quadruples of Theorem \ref{thm:Pythagorean}, namely
\begin{equation}
(p, q, r, s) = \left(\frac{5}{18}\,\pi, \frac{2}{9}\,\pi, \frac{13}{18}\,\pi, \frac{11}{18}\,\pi \right),
\end{equation}
which corresponds to item $11$ in Appendix B with parameter $t = \frac{\pi}{18}$.

The corresponding $\mathbb{Z}_2$-symmetric rational tetrahedron $T$ has edge lengths 
\begin{equation}
(\ell_p, \ell_q, \ell_r, \ell_s) = \left(\frac{5}{18}\,\pi, \frac{2}{9}\,\pi, \frac{5}{18}\,\pi, \frac{7}{18}\,\pi\right),
\end{equation}
and volume $\mathrm{vol}\, T = \pi^2/162$.

We shall prove that $T$ cannot be decomposed into any finite number of spherical Coxeter tetrahedra $S_i$, $i=1,\dots,11$, c.f. Table~\ref{tabular:coxeter-tetrahedra}. 

Suppose that it were indeed the case: then the vertex links of $T$ would be decomposed into a finite number of vertex links of Coxeter tetrahedra. The latter correspond to any of the Coxeter spherical triangles $\Delta_{2,2,n}$, $n\geq 2$, $\Delta_{2,3,3}$, $\Delta_{2,3,4}$ or $\Delta_{2,3,5}$. 

Let us consider one of the vertices $v$ of $T$ whose link $\mathrm{Lk}_v$ is a spherical triangle $\tau$ with angles $\alpha = \frac{5\pi}{18}$, $\beta = \frac{2\pi}{9}$ and $\gamma = \frac{11\pi}{18}$. The side lengths of this triangle opposite to the above mentioned angles are denoted by $\ell_\alpha$, $\ell_\beta$ and $\ell_\gamma$, respectively. The spherical law of cosines \cite[Theorem 2.5.3]{Ratcliffe} grants that $\frac{\pi}{6} < \ell_\alpha, \ell_\beta, \ell_\gamma < \frac{\pi}{2}$. We can thus position $\tau$ on the sphere $\mathbb{S}^2 = \{ (x,y,z)\in \mathbb{R}^2 | x^2 + y^2 + z^2 = 1 \}$ so that one of its vertices has coordinates $(1,0,0)$, and its adjacent vertex has coordinates $(\cos \ell_\gamma, \sin \ell_\gamma, 0)$, while the third one is in the intersection of the positive orthant $\{ (x,y,z)\in \mathbb{R}^2 | x,\, y,\, z \geq 0 \}$ with $\mathbb{S}^2$. Then we can verify that all the vertices of $\tau$ lie in the circle of radius $\frac{\pi}{4}$ centred at $p = \left(\cos\frac{4 \pi}{25}, \sin\frac{4 \pi}{25}, 0\right)$, c.f. Monty \cite{Monty}.

Thus, $\mathrm{diam}\,\mathrm{Lk}_v < \frac{\pi}{2}$, and none of the triangles $\Delta_{2,2,n}$ is a part of the decomposition of $\mathrm{Lk}_v$. The remaining cases are limited to a decomposition into $k\geq 0$ triangles of type $\Delta_{2,3,3}$, $l\geq 0$ triangles of type $\Delta_{2,3,4}$, and $m\geq 0$ triangles of type $\Delta_{2,3,5}$. Then the obvious sum of areas equality holds:
\begin{equation*}
k\, \mathrm{Area}\,\Delta_{2,3,3} + l\, \mathrm{Area}\,\Delta_{2,3,4} + m\, \mathrm{Area}\, \Delta_{2,3,5} = \mathrm{Area}\, \mathrm{Lk}_v,
\end{equation*}
which can be simplified down to
\begin{equation*}
10 k + 5 l + 2 m = \frac{20}{3}
\end{equation*}
by using the angle excess formula for the area of a spherical triangle \cite[Theorem 2.5.5]{Ratcliffe}. The latter never holds with $k,l,m \in \mathbb{Z}$. 

Another spherical rational tetrahedron $T^\prime$ with volume $\pi^2/162$ is given by the Coxeter diagram 
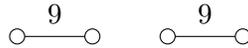
\begin{figure}[h]
\center
\begin{tikzpicture}
\draw (0,0) -- (1,0);
\draw (2,0) -- (3,0);
	
\draw[fill=white] (0,0) circle(.1);
\draw[fill=white] (1,0) circle(.1);
\draw[fill=white] (2,0) circle(.1);
\draw[fill=white] (3,0) circle(.1);
	
\node at (0.5, 0.3) {$9$};
\node at (2.5, 0.3) {$9$};
\end{tikzpicture}
\caption{The Coxeter tetrahedron $T^\prime$}\label{fig:Coxeter-tetrahedron}
\end{figure}

Both $T$ and $T^\prime$ have equal volumes and equal Dehn invariants: the former is by construction, and the latter follows from the fact that their dihedral angles are rational multiples of $\pi$, which implies that their Dehn invariants vanish. 

\begin{quest}
Are the tetrahedra $T$ and $T^\prime$, as above, scissors congruent?
\end{quest}

\section{Rational Lambert cubes}

\begin{figure}
\center
\includegraphics[scale=0.5]{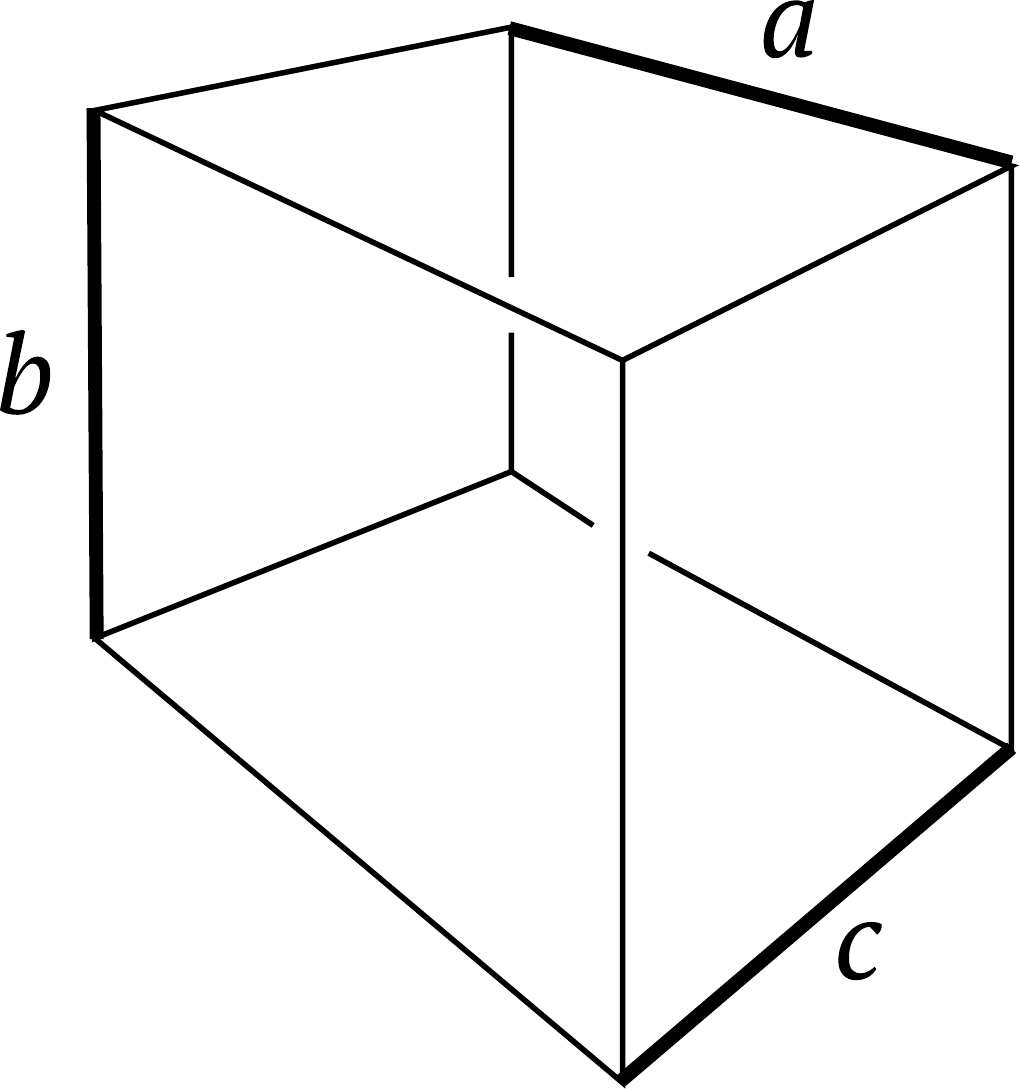}
\caption{The Lambert cube $L(a, b, c)$ with essential angles marked}\label{fig:Lambert-cube}
\end{figure}

A Lambert cube $L := L(a, b, c)$ is depicted in Figure~\ref{fig:Lambert-cube}. It is realisable as a spherical polytope $L \subset \mathbb{S}^3$, if $ \pi/2 < \alpha, \beta, \gamma < \pi$, c.f. \cite{Diaz}. All other dihedral angles of $L$, apart from the \textit{essential} ones $a$, $b$ and $c$, are always equal to $\pi/2$. 

The following fact stated as Proposition 4 in \cite{DM} holds for the volume function $\mathrm{Vol}\, L$, which allows us to seek rational Lambert cubes i.e. $L = L(a, b, c)$ with $a, b, c \in \pi\, \mathbb{Q}$, having rational volume $\mathrm{Vol}\, L \in \pi^2\, \mathbb{Q}$.

\begin{prop}[Proposition 4 in \cite{DM}]\label{prop:volume-Lambert-cube}
Suppose that the essential angles of a spherical Lambert cube $L = L(a, b, c)$ satisfy the relation $\cos^2 a + \cos^2 b + \cos^2 c = 1$. Then
\begin{equation*}
\mathrm{Vol}\, L = \frac{1}{4}\, \left(\frac{\pi^2}{2} - (\pi-a)^2 - (\pi-b)^2 - (\pi-c)^2\right).
\end{equation*}
\end{prop}

By using Monty \cite{Monty} we find, in a way analogous to the discussion in Sections~\ref{section:sporadic} -- \ref{section:families}, that there are only two sporadic rational Lambert cubes satisfying the conditions of Proposition~\ref{prop:volume-Lambert-cube}. No continuous families are present in this case, as follows from Theorem~\ref{thm:cj}.

Namely, only the following two Lambert cubes come out of our analysis: $L_1 = L(\frac{3\pi}{4}, \frac{2\pi}{3}, \frac{2\pi}{3})$ and $L_2 = L(\frac{2\pi}{3}, \frac{3\pi}{5}, \frac{4\pi}{5})$. By applying Proposition~\ref{prop:volume-Lambert-cube}, we obtain that $\mathrm{Vol}\, L_1 = 31/576\, \pi^2$ and $\mathrm{Vol}\, L_2 = 17/360\, \pi^2$. 

It is easy to produce a pair of spherical rational simplices $T_1$ and $T_2$ such that the respective $L_i$ and $T_i$, $i=1, 2$, have equal volumes and equal Dehn invariants. Let $T_1$ be given by the quadruple $(\frac{\pi}{2}, \frac{\pi}{2}, \frac{\pi}{2}, \frac{31\pi}{144})$, and let $T_2$ be given by $( \frac{\pi}{2},  \frac{\pi}{2},  \frac{\pi}{2}, \frac{17\pi}{90})$. Both $T_i$'s belong to the family $I_2(k)\times A_1^{\times 2}$ in Table~\ref{tabular:coxeter-tetrahedra}, if we allow $k$ to take rational values. 

\begin{quest}
Are the tetrahedron $T_1$ (resp. $T_2$) and the cube $L_1$ (resp. $L_2$), as above, scissors congruent?
\end{quest} 

By \cite{Felikson2}, we have that $L_1$ is the only spherical Lambert cube that can be represented as a union of mutually isometric Coxeter tetrahedra. 

\begin{quest}
Is the Lambert cube $L_2$ decomposable into any finite number of Coxeter tetrahedra?
\end{quest}

\section{Higher-dimensional aspects}

As in the proof of Theorem~\ref{thm:rational-not-coxeter-tetrahedron}, suppose that a rational $n$-dimensional, $n\geq 3$, spherical simplex $T \subset \mathbb{S}^n$ is given. Then the fact that $T$ splits into a finite number of Coxeter simplices (identified facet to facet in order to form the initial simplex $T$) will imply that all the vertex links $\mathrm{Lk}_{v_i}$, $i = 1, \dots, n+1$, of $T$ can be decomposed into a finite number of co-dimension one Coxeter simplices $T^i_j$, $j = 1, \dots, n_i$. If one of the vertex links in $T$ does \textit{not} have this property, then neither has $T$.

Let us now suppose that the three-dimensional rational tetrahedron $T^{(3)}_1 = T$ from Theorem~\ref{thm:rational-not-coxeter-tetrahedron} is a vertex link of a four-dimensional rational spherical simplex $T^{(4)}_1 \subset \mathbb{S}^4$. Then we obviously have an example of a four-dimensional simplex that does not split into any finite number of Coxeter simplices. More generally, if a rational simplex $T^{(n)}_1 \subset \mathbb{S}^n$ that is \textit{not} decomposable into Coxeter pieces can be realised as a vertex link of a rational simplex $T^{(n+1)}_1 \subset \mathbb{S}^{(n+1)}$, then $T^{(n+1)}_1$ gives us a rational simplex with the analogous property in a higher dimension. 

Constructing such a family of rational spherical simplices $T^{(n)}_1$, $n\geq 3$, starting from $T^{(3)}_1$ is simple: let $G_3 = G(T^{(3)}_1)$ be the Gram matrix of $T^{(3)}_1$, and then let $T^{(n)}_1$, $n\geq 3$, be the spherical simplex with the block-diagonal Gram matrix $$G_n := \left( \begin{array}{cccc}
G_3& \mathbf{0}& \mathbf{0}& \mathbf{0}\\
\mathbf{0}& 1& 0& 0\\
\mathbf{0}& 0& \ddots& 0\\
\mathbf{0}& 0& 0& 1
\end{array} \right).$$

The volume of $T^{(n)}_1$, $n\geq 4$, equals $\mathrm{Vol}\, T^{(n)}_1 = \mathrm{Vol}\, T^{(3)}_1 \cdot \frac{\mathrm{Vol}\, \mathbb{S}^{n}}{2^{n-3}\, \mathrm{Vol}\, \mathbb{S}^3}$, which is a rational multiple of $\mathrm{Vol}\,\mathbb{S}^n$ once $T^{(3)}_1$ has rational volume. 

If we apply the above construction to the tetrahedron $T^{(3)}_2 = T^{\prime}$, then we obtain a family of Coxeter tetrahedra $T^{(n)}_2$, each generating the respective finite reflection group $I_2(9)\times I_2(9)\times (A_1)^{n-3}$, for $n\geq 3$. The volumes and Dehn invariants of each pair $T^{(n)}_1$ and $T^{(n)}_2$, $n\geq 3$, are equal, although the former is not decomposable into any finite number of Coxeter tetrahedra, and the latter is a Coxeter tetrahedron itself. 

\begin{quest}
Are $T^{(n)}_1$ and $T^{(n)}_2$, $n \geq 3$, scissors congruent?
\end{quest} 

It is also worth mentioning that if there exists a tetrahedron $T\subset \mathbb{S}^3$ with ``rational'' dihedral angles, but ``irrational'' volume, i.e. a counterexample to the initial conjecture by Cheeger and Simons \cite{CS}, then using the above construction we can also produce a counterexample in every dimension $n \geq 3$.

\vspace*{-0.125in}

{\small
\begin{table}[!hb]
\centering
\begin{tabular}{l@{\hskip 1.5in}l}
\begin{tabular}[c]{@{}l@{}}\it Alexander Kolpakov\\ \it Institut de Math\'{e}matiques\\ \it Universit\'{e} de Neuch\^{a}tel\\ \it Suisse/Switzerland\\ \it kolpakov.alexander(at)gmail.com\end{tabular} &
\begin{tabular}[c]{@{}l@{}} \it Sinai Robins\\ \it Departamento de ci\^{e}ncia da computa\c{c}\~{a}o\\ \it Instituto de Matem\'{a}tica e Estatistica\\ \it Universidade de S\~{a}o Paulo\\ \it Brasil/Brazil \\ \it sinai.robins(at)gmail.com\end{tabular} 
\end{tabular}
\end{table}
}

\newpage

\section{Appendix A}
{\small
\begin{table}[h]
\centering
\caption{Sporadic spherical $\mathbb{Z}_2$-symmetric tetrahedra: dihedral angles have the form $(p \pi, q \pi, r \pi, s \pi)$, side lengths have the form $(\ell_p \pi, \ell_q \pi, \ell_r \pi, \ell_s \pi)$, and volumes are $v \pi^2$, with $p, q, r, s, \ell_p, \ell_q, \ell_r, \ell_s, v \in \mathbb{Q}$}
\label{table:A1}
\vspace*{0.5in}
\begin{tabular}{cccc}
no. & $(p, q, r, s)$ & $(\ell_p, \ell_q, \ell_r, \ell_s)$ & $\mathrm{Vol}$ \\
1 &  (2/3, 1/3, 3/5, 1/5) &  (2/3, 1/3, 2/5, 4/5) &  7/90 \\
2 &  (25/42, 11/42, 4/7, 2/7) &  (25/42, 11/42, 3/7, 5/7) &  67/1764 \\
3 &  (2/5, 4/15, 3/5, 8/15) &  (2/5, 4/15, 2/5, 7/15) &  19/900 \\
4 &  (2/5, 1/5, 2/3, 1/2) &  (2/5, 1/5, 1/3, 1/2) &  7/720 \\
5 &  (6/7, 2/7, 1/3, 2/7) &  (6/7, 2/7, 2/3, 5/7) &  299/1764 \\
6 &  (19/30, 17/30, 11/15, 1/3) &  (19/30, 17/30, 4/15, 2/3) &  209/900 \\
7 &  (2/3, 2/3, 4/5, 2/5) &  (2/3, 2/3, 1/5, 3/5) &  31/90 \\
8 &  (6/7, 5/7, 5/7, 2/3) &  (6/7, 5/7, 2/7, 1/3) &  1013/1764 \\
9 &  (13/30, 11/30, 11/15, 1/3) &  (13/30, 11/30, 4/15, 2/3) &  29/900 \\
10 &  (7/20, 3/20, 2/3, 3/5) &  (7/20, 3/20, 1/3, 2/5) &  17/3600 \\
11 &  (4/5, 3/5, 2/3, 1/2) &  (4/5, 3/5, 1/3, 1/2) &  59/144 \\
12 &  (23/30, 11/30, 7/15, 1/3) &  (23/30, 11/30, 8/15, 2/3) &  161/900 \\
13 &  (5/7, 1/7, 1/3, 2/7) &  (5/7, 1/7, 2/3, 5/7) &  47/1764 \\
14 &  (17/30, 11/30, 2/3, 4/15) &  (17/30, 11/30, 1/3, 11/15) &  59/900 \\
15 &  (2/3, 1/5, 2/5, 1/3) &  (2/3, 1/5, 3/5, 2/3) &  37/900 \\
16 &  (13/30, 7/30, 3/5, 1/2) &  (13/30, 7/30, 2/5, 1/2) &  67/3600 \\
17 &  (5/7, 3/7, 4/7, 1/3) &  (5/7, 3/7, 3/7, 2/3) &  335/1764 \\
18 &  (1/5, 2/15, 4/5, 11/15) &  (1/5, 2/15, 1/5, 4/15) &  1/900 \\
19 &  (31/42, 25/42, 5/7, 3/7) &  (31/42, 25/42, 2/7, 4/7) & 613/1764 \\
20 &  (11/15, 3/5, 3/5, 8/15) &  (11/15, 3/5, 2/5, 7/15) &  319/900 \\
21 &  (23/30, 13/30, 1/2, 2/5) &  (23/30, 13/30, 1/2, 3/5) & 847/3600 \\
22 &  (17/42, 11/42, 5/7, 3/7) &  (17/42, 11/42, 2/7, 4/7) & 25/1764 \\
23 &  (17/30, 7/30, 1/2, 2/5) &  (17/30, 7/30, 1/2, 3/5) &  127/3600 \\
24 &  (23/30, 19/30, 2/3, 8/15) &  (23/30, 19/30, 1/3, 7/15) &  371/900 \\
25 &  (1/3, 1/3, 4/5, 2/5) &  (1/3, 1/3, 1/5, 3/5) &  1/90 \\
26 &  (4/7, 2/7, 4/7, 1/3) &  (4/7, 2/7, 3/7, 2/3) & 83/1764 \\
27 &  (3/5, 3/5, 2/3, 2/5) &  (3/5, 3/5, 1/3, 3/5) &  109/450 \\
28 &  (1/3, 1/5, 2/3, 3/5) &  (1/3, 1/5, 1/3, 2/5) &  7/900 \\
29 &  (11/30, 7/30, 2/3, 8/15) &  (11/30, 7/30, 1/3, 7/15) & 11/900 \\
30 &  (3/5, 2/5, 3/5, 1/3) &  (3/5, 2/5, 2/5, 2/3) &  49/450
\end{tabular}
\end{table}

\begin{table}[h]
\centering
\caption{Sporadic spherical $\mathbb{Z}_2$-symmetric tetrahedra (cont.)}
\label{table:A2}
\vspace*{0.5in}
\begin{tabular}{cccc}
31 &  (13/15, 4/5, 4/5, 11/15) &  (13/15, 4/5, 1/5, 4/15) &  601/900 \\
32 &  (5/7, 4/7, 2/3, 3/7) &  (5/7, 4/7, 1/3, 4/7) &  545/1764 \\
33 &  (3/5, 4/15, 7/15, 2/5) &  (3/5, 4/15, 8/15, 3/5) &  49/900 \\
34 &  (23/30, 17/30, 3/5, 1/2) &  (23/30, 17/30, 2/5, 1/2) & 1267/3600 \\
35 &  (2/5, 2/5, 2/3, 2/5) &  (2/5, 2/5, 1/3, 3/5) &  19/450 \\
36 &  (17/20, 7/20, 2/5, 1/3) &  (17/20, 7/20, 3/5, 2/3) &  797/3600 \\
37 &  (4/5, 2/5, 1/2, 1/3) &  (4/5, 2/5, 1/2, 2/3) &  163/720 \\
38 &  (3/7, 2/7, 2/3, 3/7) &  (3/7, 2/7, 1/3, 4/7) &  41/1764 \\
39 &  (13/15, 1/5, 4/15, 1/5) &  (13/15, 1/5, 11/15, 4/5) &  91/900 \\
40 &  (19/30, 7/30, 7/15, 1/3) &  (19/30, 7/30, 8/15, 2/3) & 41/900 \\
41 &  (2/3, 2/5, 2/3, 1/5) &  (2/3, 2/5, 1/3, 4/5) &  103/900 \\
42 &  (3/5, 1/5, 1/2, 1/3) &  (3/5, 1/5, 1/2, 2/3) &  19/720 \\
43 &  (3/5, 1/3, 2/3, 1/5) &  (3/5, 1/3, 1/3, 4/5) &  43/900 \\
44 &  (4/5, 1/3, 2/5, 1/3) &  (4/5, 1/3, 3/5, 2/3) &  157/900 \\
45 &  (19/30, 13/30, 2/3, 4/15) &  (19/30, 13/30, 1/3, 11/15) &  119/900 \\
46 &  (4/5, 1/5, 1/3, 1/5) &  (4/5, 1/5, 2/3, 4/5) &  31/450 \\
47 &  (17/20, 13/20, 2/3, 3/5) &  (17/20, 13/20, 1/3, 2/5) & 1817/3600 \\
48 &  (4/5, 2/3, 4/5, 1/2) &  (4/5, 2/3, 1/5, 1/2) &  1691/3600 \\
49 &  (4/5, 2/15, 4/15, 1/5) &  (4/5, 2/15, 11/15, 4/5) &  31/900 \\
50 &  (2/5, 1/3, 4/5, 1/3) &  (2/5, 1/3, 1/5, 2/3) &  13/900 \\
51 &  (1/5, 1/5, 4/5, 2/3) &  (1/5, 1/5, 1/5, 1/3) &  1/450 \\
52 &  (2/7, 1/7, 5/7, 2/3) &  (2/7, 1/7, 2/7, 1/3) &  5/1764 \\
53 &  (11/15, 2/5, 7/15, 2/5) &  (11/15, 2/5, 8/15, 3/5) &  169/900 \\
54 &  (31/42, 17/42, 4/7, 2/7) &  (31/42, 17/42, 3/7, 5/7) & 319/1764 \\
55 &  (4/5, 4/5, 4/5, 2/3) &  (4/5, 4/5, 1/5, 1/3) &  271/450 \\
56 &  (4/5, 2/3, 2/3, 3/5) &  (4/5, 2/3, 1/3, 2/5) &  427/900 \\
57 &  (11/15, 2/3, 11/15, 1/2) &  (11/15, 2/3, 4/15, 1/2) &  493/1200 \\
58 &  (2/3, 3/5, 4/5, 1/3) &  (2/3, 3/5, 1/5, 2/3) &  253/900 \\
59 &  (13/20, 3/20, 2/5, 1/3) &  (13/20, 3/20, 3/5, 2/3) &  77/3600
\end{tabular}
\end{table}

\begin{landscape}
\section{Appendix B}

\begin{table}[h]
\centering
\caption{Continuous families of $\mathbb{Z}_2$-symmetric spherical tetrahedra}
\label{table:B1}
\begin{tabular}{ccccc}
no. & $(p, q, r, s)$ & $(\ell_p, \ell_q, \ell_r, \ell_s)$ & domain & $\mathrm{Vol}$ \\
1 &   $(1/2  \pi + t, 1/2  \pi, 1/2  \pi, 1/2  \pi )$ &   $(t +  \pi/2,  \pi/2,
 \pi/2,  \pi/2 )$ & $ $  &$(2 t +  \pi )^2/8$ \\
2 &   $(3/4  \pi - 1/2 t, 1/4  \pi - 1/2 t, 1/3  \pi - t, 1/3  \pi + t )$ & 
 $(-t/2 + 3  \pi/4, -t/2 +  \pi/4, -t + 2  \pi/3, t + 2  \pi/3 )$ &  $ $ &$-t^2/4 -
 \pi t/2 + 13  \pi^2/144$ \\
3 &   $(1/2  \pi + t, 1/2  \pi, 1/3  \pi + t, 2/3  \pi - t )$ &   $(t +  \pi/2,
 \pi/2, t +  \pi/3, -t + 2  \pi/3 )$ &  $ $ &$\pi (6 t +  \pi )/9$ \\
4 &   $(1/2  \pi, 1/6  \pi + t, 2/3  \pi - t, 1/3  \pi + t )$ &   $( \pi/2, t +
 \pi/6, -t + 2  \pi/3, t +  \pi/3 )$ &$ $ &  $ \pi t/3$ \\
5 &   $(2/3  \pi - t, 1/3  \pi, 1/3  \pi + t, 1/2  \pi )$ &   $(-t + 2  \pi/3,
 \pi/3,  \pi/2, -t + 2  \pi/3 )$ &  $ $ &$t^2/4 -  \pi t/3 + 5  \pi^2/48$ \\
6 &   $(1/2  \pi, 1/2  \pi - t, 1/3  \pi + t, 2/3  \pi - t )$ &   $( \pi/2, -t +
 \pi/2, t +  \pi/3, -t + 2  \pi/3 )$ &  $ $ &$ \pi (-3 t +  \pi )/9$ \\
7 &   $(1/3  \pi + t, 1/3  \pi, 1/2  \pi, 2/3  \pi - t )$ &   $(t +  \pi/3,
 \pi/3, t +  \pi/3,  \pi/2 )$ &  $ $ &$t^2/4 +  \pi t/6 +  \pi^2/48$ \\
8 &   $(2/3  \pi, 1/3  \pi - t, 1/2  \pi, 1/3  \pi - t )$ &   $(2  \pi/3, -t +
 \pi/3, t + 2  \pi/3,  \pi/2 )$ &  $ $ &$t^2/4 - 2  \pi t/3 + 5  \pi^2/48$ \\
9 &   $(2/3  \pi, 1/3  \pi + t, 1/3  \pi + t, 1/2  \pi )$ &   $(2  \pi/3, t +
 \pi/3,  \pi/2, -t + 2  \pi/3 )$ &  $ $ &$t^2/4 + 2  \pi t/3 + 5  \pi^2/48$ \\
10 &   $(1/2  \pi, 1/2  \pi - t, 1/3  \pi - t, 2/3  \pi + t )$ &   $( \pi/2, -t
+  \pi/2, -t +  \pi/3, t + 2  \pi/3 )$ &  $ $ &$ \pi  (-6 t +  \pi )/9$ \\
11 &   $(1/4  \pi + 1/2 t, 1/4  \pi - 1/2 t, 2/3  \pi + t, 2/3  \pi - t )$ &
 $(t/2 +  \pi/4, -t/2 +  \pi/4, t +  \pi/3, -t +  \pi/3 )$ &  $ $ &$-t^2/4 +
 \pi^2/144$ \\
12 &   $(1/2  \pi + t, 1/2  \pi, 1/3  \pi - t, 2/3  \pi + t )$ &   $(t +  \pi/2,
 \pi/2, -t +  \pi/3, t + 2  \pi/3 )$ &  $0 \leq t \leq \pi/6$ &$ \pi (3 t +  \pi )/9$ \\
13 &   $(1/2  \pi, 1/6  \pi + t, 1/2  \pi, 1/2  \pi )$ &   $( \pi/2, t +  \pi/6,
 \pi/2,  \pi/2 )$ &  for no. 1 -- no. 34 & $(6 t +  \pi )^2/72$ \\
14 &   $(1/2  \pi, 1/2  \pi - t, 1/2  \pi, 1/2  \pi )$ &   $( \pi/2, -t +  \pi/2,
 \pi/2,  \pi/2 )$ &  $ $ & $(2 t -  \pi )^2/8$ \\
15 &   $(1/3  \pi, 1/3  \pi - t, 2/3  \pi + t, 1/2  \pi )$ &   $( \pi/3, -t +
 \pi/3,  \pi/2, -t +  \pi/3 )$ &  $ $ &$t^2/4 -  \pi t/6 +  \pi^2/48$ \\
16 &   $(3/4  \pi + 1/2 t, 1/4  \pi + 1/2 t, 1/3  \pi - t, 1/3  \pi + t )$ &
 $(t/2 + 3  \pi/4, t/2 +  \pi/4, -t + 2  \pi/3, t + 2  \pi/3 )$ &  $ $ &$-t^2/4 +
 \pi t/2 + 13  \pi^2/144$ \\
17 &   $(3/4  \pi + 1/2 t, 1/4  \pi + 1/2 t, 1/3  \pi + t, 1/3  \pi - t )$ &
 $(t/2 + 3  \pi/4, t/2 +  \pi/4, t + 2  \pi/3, -t + 2  \pi/3 )$ &  $ $ &$-t^2/4 +
 \pi t/2 + 13  \pi^2/144$ \\
18 &   $(1/4  \pi + 1/2 t, 1/4  \pi - 1/2 t, 2/3  \pi - t, 2/3  \pi + t )$ &
 $(t/2 +  \pi/4, -t/2 +  \pi/4, -t +  \pi/3, t +  \pi/3 )$ &  $ $ &$-t^2/4 +
 \pi^2/144$ \\
19 &   $(2/3  \pi + t, 1/3  \pi, 1/3  \pi - t, 1/2  \pi )$ &   $(t + 2  \pi/3,
 \pi/3,  \pi/2, t + 2  \pi/3 )$ &  $ $ &$t^2/4 +  \pi t/3 + 5  \pi^2/48$ \\
20 &   $(1/2  \pi, 1/6  \pi - t, 1/2  \pi, 1/2  \pi )$ &   $( \pi/2, -t +  \pi/6,
 \pi/2,  \pi/2 )$ &   $ $ &$(6 t -  \pi )^2/72$ \\
21 &   $(2/3  \pi, 1/3  \pi + t, 1/2  \pi, 1/3  \pi + t )$ &   $(2  \pi/3, t +
 \pi/3, -t + 2  \pi/3,  \pi/2 )$ &  $ $ &$t^2/4 + 2  \pi t/3 + 5  \pi^2/48$ \\
22 &   $(3/4  \pi - 1/2 t, 1/4  \pi - 1/2 t, 1/3  \pi + t, 1/3  \pi - t )$ &
 $(-t/2 + 3  \pi/4, -t/2 +  \pi/4, t + 2  \pi/3, -t + 2  \pi/3 )$ & $ $ &$-t^2/4 -
 \pi t/2 + 13  \pi^2/144$ \\
23 &   $(1/2  \pi, 1/2  \pi - t, 2/3  \pi - t, 1/3  \pi + t )$ &   $( \pi/2, -t
+  \pi/2, -t + 2  \pi/3, t +  \pi/3 )$ &  $ $ &$ \pi (-3 t +  \pi )/9$ \\
24 &   $(1/2  \pi, 1/6  \pi + t, 1/3  \pi + t, 2/3  \pi - t )$ &   $( \pi/2, t +
 \pi/6, t +  \pi/3, -t + 2  \pi/3 )$ &  $ $ &$ \pi t/3$ \\
25 &   $(1/2  \pi + t, 1/2  \pi, 2/3  \pi + t, 1/3  \pi - t )$ &   $(t +  \pi/2,
 \pi/2, t + 2  \pi/3, -t +  \pi/3 )$ &  $ $ &$ \pi (3 t +  \pi )/9$ \\
\end{tabular}
\end{table}

\clearpage 
 
 \begin{table}[h]
\centering

\caption{Continuous families of $\mathbb{Z}_2$-symmetric spherical tetrahedra (cont.)}

\label{table:B2}
\begin{tabular}{ccccc}
26 &   $(3/4  \pi + 1/2 t, 3/4  \pi - 1/2 t, 2/3  \pi - t, 2/3  \pi + t )$ &
 $(t/2 + 3  \pi/4, -t/2 + 3  \pi/4, -t +  \pi/3, t +  \pi/3 )$ &  $ $ &$-t^2/4 +
73  \pi^2/144$ \\
27 &   $(1/3  \pi + t, 1/3  \pi, 2/3  \pi - t, 1/2  \pi )$ &   $(t +  \pi/3,
 \pi/3,  \pi/2, t +  \pi/3 )$ &  $ $ &$t^2/4 +  \pi t/6 +  \pi^2/48$ \\
28 &   $(2/3  \pi - t, 1/3  \pi, 1/2  \pi, 1/3  \pi + t )$ &   $(-t + 2  \pi/3,
 \pi/3, -t + 2  \pi/3,  \pi/2 )$ &  $ $ &$t^2/4 -  \pi t/3 + 5  \pi^2/48$ \\
29 &   $(1/3  \pi, 1/3  \pi - t, 1/2  \pi, 2/3  \pi + t )$ &   $( \pi/3, -t +
 \pi/3, -t +  \pi/3,  \pi/2 )$ &  $ $ &$t^2/4 -  \pi t/6 +  \pi^2/48$ \\
30 &   $(1/2  \pi, 1/2  \pi - t, 2/3  \pi + t, 1/3  \pi - t )$ &   $( \pi/2, -t
+  \pi/2, t + 2  \pi/3, -t +  \pi/3 )$ &  $ $ &$ \pi (-6 t +  \pi )/9$ \\
31 &   $(2/3  \pi + t, 1/3  \pi, 1/2  \pi, 1/3  \pi - t )$ &   $(t + 2  \pi/3,
 \pi/3, t + 2  \pi/3,  \pi/2 )$ &  $ $ &$t^2/4 +  \pi t/3 + 5  \pi^2/48$ \\
32 &   $(2/3  \pi, 1/3  \pi - t, 1/3  \pi - t, 1/2  \pi )$ &   $(2  \pi/3, -t +
 \pi/3,  \pi/2, t + 2  \pi/3 )$ &  $ $ &$t^2/4 - 2  \pi t/3 + 5  \pi^2/48$ \\
33 &   $(1/2  \pi + t, 1/2  \pi, 2/3  \pi - t, 1/3  \pi + t )$ &   $(t +  \pi/2,
 \pi/2, -t + 2  \pi/3, t +  \pi/3 )$ &  $ $ &$ \pi (6 t +  \pi )/9$ \\
34 &   $(3/4  \pi + 1/2 t, 3/4  \pi - 1/2 t, 2/3  \pi + t, 2/3  \pi - t )$ &
 $(t/2 + 3  \pi/4, -t/2 + 3  \pi/4, t +  \pi/3, -t +  \pi/3 )$ &  $ $ &$-t^2/4 +
73  \pi^2/144$ \\
35 &  $(1/2 \pi, 1/2 \pi - u,  \pi - t, t)$ &  $( \pi/2, -u +  \pi/2, -t +
 \pi, t)$ &  Domain A &$-t^2/2 +  \pi t/2 + u^2/2 -  \pi u/2$ \\
36 &  $(1/2  \pi, 1/2  \pi - u, t,  \pi - t)$ &  $( \pi/2, -u +  \pi/2, t, -t
+  \pi)$ & Domain A &$-t^2/2 +  \pi t/2 + u^2/2 -  \pi u/2$ \\
37 &  $(1/2 \pi + u, 1/2 \pi,  \pi - t, t)$ &  $(u +  \pi/2,  \pi/2, -t +
 \pi, t$) &  Domain A &$-t^2/2 +  \pi t/2 + u^2/2 +  \pi u/2$ \\
38 &  $(1/2  \pi + u, 1/2  \pi, t,  \pi - t)$ &  $(u +  \pi/2,  \pi/2, t, -t +
 \pi)$ & Domain A &$ -t^2/2 +  \pi t/2 + u^2/2 +  \pi u/2$ \\
39 &  $(1/2  \pi, 1/2  \pi - t,  \pi - u, u)$ &  $( \pi/2, -t +  \pi/2, -u +
 \pi, u)$ & Domain B &$ t^2/2 -  \pi t/2 - u^2/2 +  \pi u/2$ \\
40 &  $(1/2  \pi, 1/2  \pi - t, u,  \pi - u)$ &  $( \pi/2, -t +  \pi/2, u, -u
+  \pi)$ & Domain B &$ t^2/2 -  \pi t/2 - u^2/2 +  \pi u/2$ \\
41 &  $(1/2  \pi + t, 1/2  \pi,  \pi - u, u)$ &  $(t +  \pi/2,  \pi/2, -u +
 \pi, u)$ &  Domain B &$t^2/2 +  \pi t/2 - u^2/2 +  \pi u/2$ \\
42 &  $(1/2  \pi + t, 1/2  \pi, u,  \pi - u)$ &  $(t +  \pi/2,  \pi/2, u, -u +
 \pi)$ &  Domain B &$t^2/2 +  \pi t/2 - u^2/2 +  \pi u/2$ \\
\end{tabular}
\end{table}

\centering
\begin{tabular}{c|c}
Domain A: & Domain B: \\
$0 \leq u \leq \frac{\pi}{2}$, & $0 \leq u \leq \pi$, \\
$0 \leq t \leq \pi$, & $0 \leq t \leq \frac{\pi}{2}$, \\
$t \geq u$. & $t \leq u$.\\ 
\end{tabular}

\end{landscape}
}
\end{document}